\documentclass[conference]{IEEEtran}
\IEEEoverridecommandlockouts

\usepackage{xcolor}
\definecolor{red}{rgb}{1,0,0}
\definecolor{orange}{rgb}{1,.5,0}
\definecolor{darkgreen}{rgb}{0,.5,0}
\definecolor{grey}{rgb}{0.5,0.5,0.5}
\newcommand{\PH}[1]{\textcolor{red}{Philipp: #1}}

\newcommand{\minn}{\mathop{\rm min}}
\newcommand{\argminn}{\mathop{\rm argmin}}
\newcommand{\maxn}{\mathop{\rm max}}
\newcommand{\reals}{{\mbox{\bf R}}}

\newcommand{\st}{\rm s.t.}
\ifCLASSINFOpdf
  \usepackage[pdftex]{graphicx}
  \DeclareGraphicsExtensions{.pdf,.jpeg,.png}
\else
\fi
\usepackage{amsmath}
\ifCLASSOPTIONcompsoc
 \usepackage[caption=false,font=normalsize,labelfont=sf,textfont=sf]{subfig}
\else
 \usepackage[caption=false,font=footnotesize]{subfig}
\fi
\usepackage{xprintlen}

\usepackage{stfloats}
\newcommand{\Hquad}{\hspace{0.5em}} 

\begin{document}
%
\title{Bilevel Optimization for Improved Flexibility Aggregation Models of Electric Vehicle Fleets}
%
%
%


\author{
\IEEEauthorblockN{Philipp Härtel}
\IEEEauthorblockA{\textit{Fraunhofer Institute for Energy Economics} \\
\textit{and Energy System Technology, Fraunhofer IEE}\\
\textit{Sustainable Electrical Energy Systems, University of Kassel}\\
Kassel, Germany \\
philipp.haertel@iee.fraunhofer.de}
\and
\IEEEauthorblockN{Michael von Bonin}
\IEEEauthorblockA{\textit{Fraunhofer Institute for Energy Economics} \\
\textit{and Energy System Technology, Fraunhofer IEE}\\
Kassel, Germany \\
michael.von.bonin@iee.fraunhofer.de\\
}
}


%
%

\markboth{Journal of \LaTeX\ Class Files,~Vol.~14, No.~8, August~2015}%
{Shell \MakeLowercase{\textit{et al.}}: A Novel Aggregation Method for Electric V}
%



\maketitle

\begin{abstract}
Electric vehicle (EV) fleets are expected to become an increasingly important source of flexibility for power system operations.
However, accurately capturing the flexibility potential of numerous and heterogeneous EVs remains a significant challenge.
We propose a bilevel optimization formulation to enhance flexibility aggregations of electric vehicle fleets.
The outer level minimizes scheduling deviations between the aggregated and reference EV units, while the inner level maximizes the aggregated unit's profits.
Our approach introduces hourly to daily scaling factor mappings to parameterize the aggregated EV units.
Compared to simple aggregation methods, the proposed framework reduces the root-mean-square error of charging power by 78~per cent, providing more accurate flexibility representations.
The proposed framework also provides a foundation for several potential extensions in future work.
\end{abstract}


%
\IEEEpeerreviewmaketitle

\section*{Nomenclature}

\subsection{Sets and indices}
\addcontentsline{toc}{section}{Nomenclature}
\begin{IEEEdescription}[\IEEEusemathlabelsep\IEEEsetlabelwidth{$V_1,V_2$}]
\item[$\mathcal{T}$, $t$] Set of discrete time steps, indexed by $t$. $\mathcal{T}=\{1,\dots,T\}$, $\mathcal{T}^\prime =\mathcal{T}\setminus \{T\}$, $\mathcal{T}^{\prime\prime} =\mathcal{T}\setminus \{1\}$.
\item[$\mathcal{V}$, $v$] Set of individual EV units, indexed by $v$.
\item[$\mathcal{U}$, $u$] Set of AEV units, indexed by $u$.
\item[$\mathcal{V}_u$, $v$] Set of EV units aggregated by AEV unit $u$.
\end{IEEEdescription}

\subsection{Parameters}
\addcontentsline{toc}{section}{Nomenclature}
\begin{IEEEdescription}[\IEEEusemathlabelsep\IEEEsetlabelwidth{$\underline{X}^{\text{C}}_{v,t}/\overline{X}^{\text{C}}_{v,t}$}]
\item[$\Pi_t$] El. price (EUR/MWh).
\item[$\Phi^\text{DR}_{v,t}$] El. demand for driving (MWh/h).
\item[$\Phi^\text{TH}_{v,t}$] El. demand for thermal mngt. (MWh/h).
\item[$\rho_{v}$] Self-discharge factor (1).
\item[$\eta_{v}^\text{C/D}$] Charge/Discharge efficiency (1).
\item[$\underline{X}^{\text{C}}_{v,t}/\overline{X}^{\text{C}}_{v,t}$] Min./Max. charging availability (MWh/h).
\item[$\underline{X}^{\text{D}}_{v,t}/\overline{X}^{\text{D}}_{v,t}$] Min./Max. discharging availability (MWh/h).
\item[$\underline{X}^{\text{S}}_{v,t}/\overline{X}^{\text{S}}_{v,t}$] Min./Max. storage level (SOC) (MWh).
\item[$\gamma^{(\cdot)}$] Deviation weights (1).
\end{IEEEdescription}

\subsection{Primal decision variables}
\addcontentsline{toc}{section}{Nomenclature}
\begin{IEEEdescription}[\IEEEusemathlabelsep\IEEEsetlabelwidth{$\overline{\kappa}^{\text{(.)}}_{u,\tau},\underline{\kappa}^{\text{(.)}}_{u,\tau}$}]
\item[$x^{\text{C}}_{v,t}$] Single EV charging power (MWh/h).
\item[$x^{\text{D}}_{v,t}$] Single EV discharging power (MWh/h).
\item[$x^{\text{S}}_{v,t}$] Single EV storage level (SOC) (MWh).
\item[$x^{\text{C}}_{u,t}$] AEV charging power (MWh/h).
\item[$x^{\text{D}}_{u,t}$] AEV discharging power (MWh/h).
\item[$x^{\text{S}}_{u,t}$] AEV storage level (SOC) (MWh).
\item[$\underline{\kappa}^{\text{(.)}}_{u,\tau},\overline{\kappa}^{\text{(.)}}_{u,\tau}$] AEV scaling factors (1).
\item[$\underline{x}^{\text{C}}_{u,t},\overline{x}^{\text{C}}_{u,t}$] Min./Max. AEV charging avail. (MWh/h).
\item[$\underline{x}^{\text{D}}_{u,t},\overline{x}^{\text{D}}_{u,t}$] Min./Max. AEV discharging avail. (MWh/h).
\item[$\underline{x}^{\text{S}}_{u,t},\overline{x}^{\text{S}}_{u,t}$] Min./Max. AEV storage level (SOC) (MWh/h).
\end{IEEEdescription}

\subsection{Dual decision variables}
\addcontentsline{toc}{section}{Nomenclature}
\begin{IEEEdescription}[\IEEEusemathlabelsep\IEEEsetlabelwidth{$\mu^{(\cdot)}_{u,t}$}]
\item[$\lambda_{u,t}$] Dual variables of AEV storage continuity equation.
\item[$\mu^{(\cdot)}_{u,t}$] Dual variables of AEV bound inequalities.
\end{IEEEdescription}

\section{Introduction}
%
%
%
%

Integrating a growing number of electric vehicles (EVs) into electricity markets and power grids presents challenges and opportunities for effective and efficient system design and operation~\cite{Braun.2024}.
EV charging imposes additional demand on the grid, requiring greater flexibility to accommodate varying charging patterns and behaviors. 
\begin{figure}[b]
\centering
\includegraphics[width=\columnwidth]{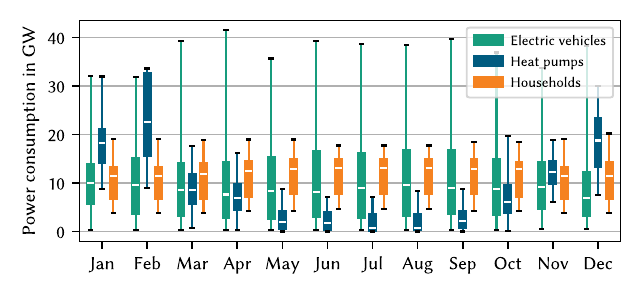}
\caption{Monthly distributions of hourly power consumption for electric vehicles, heat pumps, and household loads in Germany under a net-neutral scenario. Boxplots show the inter-quartile range, with whiskers extending to the 2.5\,\% and 97.5\,\% percentiles.} 
\label{fig:figure_boxplots}
\end{figure}

Fig.~\ref{fig:figure_boxplots} illustrates the monthly distribution of electric power demand for electric vehicles, heat pumps, and conventional household loads in Germany’s 2045 electricity market, assuming the widespread adoption of flexible consumers.
Seasonal variations are evident, with peak consumption in winter.
While the median power demand for EVs remains below 10\,GW year-round, the upper whisker often exceeds 40\,GW, showing substantial variability that needs to be handled appropriately in model aggregations.
Introducing vehicle-to-grid (V2G) capabilities further amplifies the need for accurate flexibility representations in system planning tools.
To address these challenges, we propose a bilevel optimization approach for identifying the flexibility aggregation of EVs for power and energy system planning models.

\subsection{Electric vehicle aggregation problem}
Simulating and analyzing large and highly integrated power and energy systems easily becomes computationally intractable.
Using model aggregations, projections, simplifications, or abstractions with a smaller overall problem size helps approach this intractability~\cite{Rogers.1991}.
Methods that develop standalone or auxiliary models with reduced complexity but simultaneously provide reasonable approximations of the original problem are instrumental~\cite{Hartel.2021}.
Hence, when modeling large EV fleets, aggregation plays a pivotal role.
By consolidating the granular details of individual EVs into a representative surrogate model, aggregation enables the analysis of fleet behavior at a higher level.
Finding a good trade-off between computational burden reduction and model accuracy is a key challenge since the aggregation methods influence the quality of insights and decisions derived from energy system models.
\par
Na\"ive aggregation methods, such as single-vehicle aggregation, estimate fleet flexibility by summing individual EVs' charging availability and driving demands.
While computationally simple, these approaches often oversimplify EV behavior, neglecting vehicle interactions and dependencies.
Consequently, such methods risk overestimating the flexibility of the fleet, leading to unrealistic system planning assessments.
\par
Advanced aggregation techniques have been developed to address these limitations and capture the dynamics and heterogeneity of EV fleets.
One such method is the virtual battery aggregation approach, which models EV fleets as virtual energy storage systems.
It considers charging and discharging limits, battery capacities, and individual EV characteristics (e.g., driving demands) to estimate the fleet's combined flexibility.
By incorporating the flexibility of each vehicle, this method provides a more realistic representation of aggregated EV fleets.
It serves as a reference for the bilevel optimization approach proposed in this paper.
Applications of this methodology are demonstrated in \cite{Muessel.2023}, \cite{Trost.2017}, and \cite{Pertl.2019}.
\par
Clustering techniques group EVs with similar charging patterns and characteristics to avoid mixing heterogeneous vehicles during the aggregation process.
By aggregating vehicles within each cluster, this approach improves the accuracy of fleet flexibility estimation \cite{Geng.2023}, \cite{Mukhi.2024} and can be combined with other methods to enhance overall accuracy and efficiency.
\par
While different methods have been proposed to build reasonable aggregations of EV fleet flexibility, many approaches are too coarse, relying on hand-tuned and general scaling factors, leading to flexibility descriptions that under- or overestimate the flexibility potential.

\subsection{Bilevel optimization for EV Aggregation}
We adopt an approach that has been prominent in identifying accurate reduced models for hydropower systems.
Since operating hydropower for complex river systems, various connections of plants and reservoirs, and non-convex production characteristics is computationally challenging, ``composite'' or ``equivalent'' models are commonly used to simulate its operation in power system operation models \cite{Almeida.2013}, \cite{Risberg.2017}, \cite{Blom.2020}, \cite{Schmitz.2023}.
The core idea is to use a bilevel problem formulation to compute the aggregated model.
The bilevel problem is an optimization problem constrained by another optimization problem, resulting in an outer and inner level~\cite{Dempe.2015}.
\par
\begin{figure}
\centering
\includegraphics[width=\columnwidth]{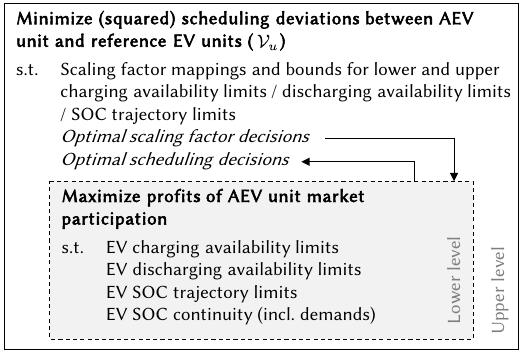}
\caption{Schematic overview of the bilevel optimization problem structure.}
\label{fig:bilevel_structure}
\end{figure}

We use a bilevel optimization approach to minimize the difference in simulated charging and discharging schedule decisions between an aggregated electric vehicle (AEV) unit and the individual reference EV units.
Fig.~\ref{fig:bilevel_structure} illustrates the underlying bilevel optimization problem structure.
The key decision variables of interest are the scaling factors to characterize the EV fleet flexibility represented by the AEV unit.
\par
The remainder of this work is organized as follows: Section II introduces the individual EV unit dispatch and the novel bilevel AEV unit problem formulations, including the necessary reformulations.
Section III presents a numerical case study to illustrate the effectiveness of the proposed approach.
Finally, Section IV summarizes the key findings, discusses their implications, and highlights future research directions in EV aggregation.

\section{Methodology}

\subsection{Key assumptions}
Our bilevel approach is based on two fundamental assumptions.
First, we assume that each individual EV interacts with the electricity market at a given price. Energy system models typically assume this implicitly (dualized) via market clearing constraints.
Second, we assume that the flexibility characterization requires that each scaling factor applies to at least two time steps with different realizations of other input parameters, e.g., price, charging availability, state-of-charge, and demands, to not just exactly match the reference profile.

\subsection{Individual EV unit dispatch problem}
For a given individual vehicle $v\in \mathcal{V}$, we write the individual EV unit dispatch problem as the following single-level optimization problem:
\begin{subequations}{\allowdisplaybreaks
\begin{IEEEeqnarray}{Crlc}
\minn_{\Xi_\text{SL}}  & \IEEEeqnarraymulticol{3}{l}{\sum_{t\in \mathcal{T}} \Big [ \Pi_t \left(x^\text{C}_{v,t} - x^\text{D}_{v,t}\right) \Big ]} \\
\st  & \underline{X}^{\text{C}}_{v,t} &\le x^{\text{C}}_{v,t} \le \overline{X}^{\text{C}}_{v,t}                   &\quad \forall t\in \mathcal{T}\label{eq:ivp_ch_bd}\\
     & \underline{X}^{\text{D}}_{v,t} &\le x^{\text{D}}_{v,t} \le \overline{X}^{\text{D}}_{v,t}                   &\quad \forall t\in \mathcal{T}\label{eq:ivp_dch_bd}\\
     & \underline{X}^\text{S}_{v,t} &\le x^\text{S}_{v,t} \le \overline{X}^\text{S}_{v,t}                         &\quad \forall t\in \mathcal{T}\label{eq:ivp_s_bd}\\
     & x^\text{S}_{v,t+1} &= \rho_{v} x^\text{S}_{v,t} - \Phi^\text{DR}_{v,t} - \Phi^\text{TH}_{v,t}                  &\nonumber\\
     & & \IEEEeqnarraymulticol{2}{r}{+ \eta_{v}^\text{C} x^{\text{C}}_{v,t} - \textstyle\frac{1}{\eta_{v}^\text{D}} x^{\text{D}}_{v,t} \quad \forall t\in \mathcal{T}^\prime}\label{eq:ivp_s_bal}\\
     & X^\text{S}_{v,T+1} &= \rho_{v} x^\text{S}_{v,T} - \Phi^\text{DR}_{v,T} - \Phi^\text{TH}_{v,T}                  &\nonumber\\
     & & + \eta_{v}^\text{C} x^{\text{C}}_{v,T} - \textstyle\frac{1}{\eta_{v}^\text{D}} x^{\text{D}}_{v,T} & \label{eq:ivp_s_final}
\end{IEEEeqnarray}}
\end{subequations}
where $\Pi_t \in \reals\; \forall t \in \mathcal{T}$ is the electricity price, $\Phi^\text{DR}_{v,t} \in \reals\; \forall t \in \mathcal{T}$ is the electricity demand for driving, $\Phi^\text{TH}_{v,t} \in \reals\; \forall t \in \mathcal{T}$ is the electricity demand for thermal management, $\rho_{v} \in [0,1]$ is the self-discharge factor, and $\eta_{v}^\text{C} \in [0,1]$ and $\eta_{v}^\text{D} \in [0,1]$ are the charge and discharge efficiencies, respectively.
The constraints include variable bounds for charging in \eqref{eq:ivp_ch_bd}, discharging in \eqref{eq:ivp_dch_bd}, and storage scheduling in \eqref{eq:ivp_s_bd}.
The storage continuity is enforced by \eqref{eq:ivp_s_bal} and meeting final storage levels is formulated in \eqref{eq:ivp_s_final}.
We denote the set of single-level (SL) (primal) variables as $\Xi_\text{SL}$.
Solving the individual dispatch problem returns the individual dispatch solutions, which serve as the reference in the bilevel approach.

\subsection{Bilevel AEV unit problem}

\subsubsection{Outer problem}\label{subsec:outerProblem}
The outer (or upper-level) problem seeks to minimize the deviation of the proxy electric vehicle model from the aggregated dispatch schedules of individual vehicles.
The set of outer-level (OL) (primal) variables are $\Xi_\text{OL}$ as well as all inner primal and dual $\Xi_\text{IL}$ variables, which will be defined in Subsection~\ref{subsec:innerProblem}.
For a given set $\mathcal{V}_u$ of individual vehicles EVs to be aggregated, the outer level problem $\mathcal{P}^\text{OL}$ writes as follows.
\begin{subequations}
{\allowdisplaybreaks
\begin{IEEEeqnarray}{CrCc}
\minn_{\Xi_\text{OL}}  & \IEEEeqnarraymulticol{3}{l}{\sum_{t\in \mathcal{T}} \bigg [ \gamma^\text{C}\left(x^\text{C}_{u,t} - \hat{X}^\text{C}_{\mathcal{V}_u,t} \right )^2 + \gamma^\text{D}\left ( x^\text{D}_{u,t}-\hat{X}^\text{D}_{\mathcal{V}_u,t}\right ) ^2  } \nonumber \\
& \IEEEeqnarraymulticol{3}{l}{\qquad + \gamma^\text{S}\left(x^\text{S}_{u,t} - \hat{X}^\text{S}_{\mathcal{V}_u,t} \right )^2 \bigg ] }\\
\st & \underline{x}^{\text{C}}_{u,t} = \underline{\kappa}^{\text{C}}_{u,\tau} \, \sum_{v\in \mathcal{V}_u} \underline{X}^{\text{C}}_{v,t}   & \Hquad \tau = f_n^\text{C}(t), \,     &\forall t\in \mathcal{T}\\
    & \overline{x}^{\text{C}}_{u,t} = \overline{\kappa}^{\text{C}}_{u,\tau} \, \sum_{v\in \mathcal{V}_u} \overline{X}^{\text{C}}_{v,t}      & \Hquad \tau = f_n^\text{C}(t), \,     &\forall t\in \mathcal{T}\\
    & \underline{x}^{\text{D}}_{u,t} = \underline{\kappa}^{\text{D}}_{u,\tau} \, \sum_{v\in \mathcal{V}_u} \underline{X}^{\text{D}}_{v,t}   & \Hquad \tau = f_n^\text{D}(t), \,     &\forall t\in \mathcal{T}\\
    & \overline{x}^{\text{D}}_{u,t} = \overline{\kappa}^{\text{D}}_{u,\tau} \, \sum_{v\in \mathcal{V}_u} \overline{X}^{\text{D}}_{v,t}      & \Hquad \tau = f_n^\text{D}(t), \,     &\forall t\in \mathcal{T}\\
    & \underline{x}^{\text{S}}_{u,t} = \underline{\kappa}^{\text{S}}_{u,\tau} \, \sum_{v\in \mathcal{V}_u} \underline{X}^{\text{S}}_{v,t}   & \Hquad \tau = f_n^\text{S}(t), \,     &\forall t\in \mathcal{T}\label{eq:socminscaling}\\
    & \overline{x}^{\text{S}}_{u,t} = \overline{\kappa}^{\text{S}}_{u,\tau} \, \sum_{v\in \mathcal{V}_u} \overline{X}^{\text{S}}_{v,t}      & \Hquad \tau = f_n^\text{S}(t), \,     &\forall t\in \mathcal{T}\\
    & \IEEEeqnarraymulticol{3}{l}{ \underline{\kappa}^{\text{C}}_{u,\tau}, \overline{\kappa}^{\text{C}}_{u,\tau}, \underline{\kappa}^{\text{D}}_{u,\tau}, \overline{\kappa}^{\text{D}}_{u,\tau}, \underline{\kappa}^{\text{S}}_{u,\tau}, \overline{\kappa}^{\text{S}}_{u,\tau} \in \reals^+}\\
    & \IEEEeqnarraymulticol{3}{l}{ x^\text{C}_{u,t}, x^\text{D}_{u,t}, x^\text{S}_{u,t} \in \argminn_{\Xi_\text{IL}} \big \{ \mathcal{P}^\text{IL} \big \}} 
\end{IEEEeqnarray}
}
\end{subequations}
Here, $\gamma^\text{C}\in \reals^+$, $\gamma^\text{D}\in \reals^+$, $\gamma^\text{S}\in \reals^+$ are the deviation weights, $f_n^\text{C}(t)$, $f_n^\text{D}(t)$, $f_n^\text{S}(t)$ are the scaling factor mapping functions for charging, discharging, and state-of-charge, respectively. For simplicity, we assume three identical mapping functions for charging, discharging, and storage levels.
\begin{equation}
    f_n^\text{C}(t) = f_n^\text{D}(t)= f_n^\text{S}(t)= \left\lfloor \frac{t \mod 7 \cdot 24}{n} \right\rfloor \label{eq:mapping}
\end{equation}
Here, the mapping ensures that the mapping cycles every week (168 hours) while grouping every $n$ hours within the cycle, implying that each weekday is treated.
Different mappings are possible, though the mapping function becomes formally more complex if some days are grouped together, e.g., workdays. 

\subsubsection{Inner problem}\label{subsec:innerProblem}
The inner (or lower-level) problem $\mathcal{P}^\text{IL}$ representing the individual dispatch problem of the AEV unit $u\in U$ is given by Eqs. (\ref{eq:innerobj}) to (\ref{eq:innerstorageend}) below.
Note that the lower and upper bounds of the charging, discharging, and storage levels become parameters determined by the outer level problem. 
We denote the set of inner-level (IL) variables as $\Xi_\text{IL}$.
\begin{subequations}{\allowdisplaybreaks
\begin{IEEEeqnarray}{Crlrl}
\minn_{\Xi_\text{IL}}  & \IEEEeqnarraymulticol{3}{l}{\sum_{t\in \mathcal{T}} \Big [ \Pi_t \left ( x^\text{C}_{u,t} - x^\text{D}_{u,t}\right ) \Big]} \label{eq:innerobj} \\
\st  & \underline{x}^{\text{C}}_{u,t} &\le x^{\text{C}}_{u,t} \le \overline{x}^{\text{C}}_{u,t}    \;\colon\; \underline{\mu}^\text{C}_{u,t},\overline{\mu}^\text{C}_{u,t}     &\forall t\in \mathcal{T}\\
     & \underline{x}^{\text{D}}_{u,t} &\le x^{\text{D}}_{u,t} \le \overline{x}^{\text{D}}_{u,t}    \;\colon\; \underline{\mu}^\text{D}_{u,t},\overline{\mu}^\text{D}_{u,t}     &\forall t\in \mathcal{T}\\
     & \underline{x}^\text{S}_{u,t} &\le x^\text{S}_{u,t} \le \overline{x}^\text{S}_{u,t}                                     \;\colon\; \underline{\mu}^\text{S}_{u,t},\overline{\mu}^\text{S}_{u,t}     &\forall t\in \mathcal{T}\\
     & x^\text{S}_{u,t+1} &= \rho_{u} x^\text{S}_{u,t} - \Phi^\text{DR}_{u,t} - \Phi^\text{TH}_{u,t}                 &\nonumber\\
     & & \IEEEeqnarraymulticol{2}{r}{+ \eta_{u}^\text{C} x^{\text{C}}_{u,t} - \textstyle\frac{1}{\eta_{u}^\text{D}} x^{\text{D}}_{u,t} \;\colon\; \lambda_{u,t} \quad \forall t\in \mathcal{T}^\prime}\\
     & X^\text{S}_{u,T+1} &= \rho_{u} x^\text{S}_{u,T} - \Phi^\text{DR}_{u,T} - \Phi^\text{TH}_{u,T}                                  &\nonumber\\
     & & + \eta_{u}^\text{C} x^{\text{C}}_{u,T} - \textstyle\frac{1}{\eta_{u}^\text{D}} x^{\text{D}}_{u,T} \;\colon\; \lambda_{u,T} \label{eq:innerstorageend} & 
\end{IEEEeqnarray}}
\end{subequations}

\subsection{Recasting the bilevel problem}
Given the convex nature of the inner problem, we recast the bilevel problem as a single-level problem by replacing the inner problem with its Karush-Kuhn-Tucker optimality conditions:
\begin{subequations}
{\allowdisplaybreaks
\begin{IEEEeqnarray}{CrCc}
\minn_{\Xi_\text{OLIL}}  & \IEEEeqnarraymulticol{3}{l}{\sum_{t\in \mathcal{T}} \bigg [ \gamma^\text{C}\left(x^\text{C}_{u,t} - \hat{X}^\text{C}_{\mathcal{V}_u,t} \right )^2 + \gamma^\text{D}\left ( x^\text{D}_{u,t}-\hat{X}^\text{D}_{\mathcal{V}_u,t}\right ) ^2  } \nonumber \\
& \IEEEeqnarraymulticol{3}{l}{\qquad + \gamma^\text{S}\left(x^\text{S}_{u,t} - \hat{X}^\text{S}_{\mathcal{V}_u,t} \right )^2 \bigg ] }\\
\st & \underline{x}^{\text{C}}_{u,t} = \underline{\kappa}^{\text{C}}_{u,\tau} \, \sum_{v\in \mathcal{V}_u} \underline{X}^{\text{C}}_{v,t}   & \Hquad \tau = f^\text{C}(t), \;     &\forall t\in \mathcal{T}\;\\
    & \overline{x}^{\text{C}}_{u,t} = \overline{\kappa}^{\text{C}}_{u,\tau} \, \sum_{v\in \mathcal{V}_u} \overline{X}^{\text{C}}_{v,t}      & \Hquad \tau = f^\text{C}(t), \;     &\forall t\in \mathcal{T}\;\\
    & \underline{x}^{\text{D}}_{u,t} = \underline{\kappa}^{\text{D}}_{u,\tau} \, \sum_{v\in \mathcal{V}_u} \underline{X}^{\text{D}}_{v,t}   & \Hquad \tau = f^\text{D}(t), \;     &\forall t\in \mathcal{T}\;\\
    & \overline{x}^{\text{D}}_{u,t} = \overline{\kappa}^{\text{D}}_{u,\tau} \, \sum_{v\in \mathcal{V}_u} \overline{X}^{\text{D}}_{v,t}      & \Hquad \tau = f^\text{D}(t), \;     &\forall t\in \mathcal{T}\;\\
    & \underline{x}^{\text{S}}_{u,t} = \underline{\kappa}^{\text{S}}_{u,\tau} \, \sum_{v\in \mathcal{V}_u} \underline{X}^{\text{S}}_{v,t}   & \Hquad \tau = f^\text{S}(t), \;     &\forall t\in \mathcal{T}\;\\
    & \overline{x}^{\text{S}}_{u,t} = \overline{\kappa}^{\text{S}}_{u,\tau} \, \sum_{v\in \mathcal{V}_u} \overline{X}^{\text{S}}_{v,t}      & \Hquad \tau = f^\text{S}(t), \;     &\forall t\in \mathcal{T}\;\\
    & \IEEEeqnarraymulticol{3}{l}{ \underline{\kappa}^{\text{C}}_{u,\tau}, \overline{\kappa}^{\text{C}}_{u,\tau}, \underline{\kappa}^{\text{D}}_{u,\tau}, \overline{\kappa}^{\text{D}}_{u,\tau}, \underline{\kappa}^{\text{S}}_{u,\tau}, \overline{\kappa}^{\text{S}}_{u,\tau} \in \reals^+}\\
    & \underline{x}^{\text{C}}_{u,t} \le x^{\text{C}}_{u,t} \le \overline{x}^{\text{C}}_{u,t}  &  &\forall t\in \mathcal{T}\\
    & \underline{x}^{\text{D}}_{u,t} \le x^{\text{D}}_{u,t} \le \overline{x}^{\text{D}}_{u,t}  &  &\forall t\in \mathcal{T}\\
    & \underline{x}^{\text{S}}_{u,t} \le x^{\text{S}}_{u,t} \le \overline{x}^{\text{S}}_{u,t}  &  &\forall t\in \mathcal{T}\\
    & \IEEEeqnarraymulticol{3}{l}{x^\text{S}_{u,t+1} = \rho_{u} x^\text{S}_{u,t} - \Phi^\text{DR}_{u,t} - \Phi^\text{TH}_{u,t}} \nonumber\\
    & \IEEEeqnarraymulticol{2}{r}{+ \eta_{u}^\text{C} x^{\text{C}}_{u,t} - \textstyle\frac{1}{\eta_{u}^\text{D}} x^{\text{D}}_{u,t}} &\forall t\in \mathcal{T}^\prime\\
    & \IEEEeqnarraymulticol{3}{l}{X^\text{S}_{u,T+1} = \rho_{u} x^\text{S}_{u,T} - \Phi^\text{DR}_{u,T} - \Phi^\text{TH}_{u,T}} \nonumber\\
    & \IEEEeqnarraymulticol{2}{r}{+ \eta_{u}^\text{C} x^{\text{C}}_{u,T} - \textstyle\frac{1}{\eta_{u}^\text{D}} x^{\text{D}}_{u,T}} & \\
    & \IEEEeqnarraymulticol{2}{l}{\Pi_t + \mu^1_{u,t} - \mu^2_{u,t} + \eta^\text{C} \lambda_{u,t} = 0}    &\forall t\in \mathcal{T}\\
    & \IEEEeqnarraymulticol{2}{l}{-\Pi_t + \mu^3_{u,t} - \mu^4_{u,t} - \textstyle\frac{1}{\eta_{u}^\text{D}} \lambda_{u,t} = 0}   &\forall t\in \mathcal{T}\\
    & \IEEEeqnarraymulticol{2}{l}{\mu^5_{u,t} - \mu^6_{u,t} + \rho_{u} \lambda_{u,t} - \lambda_{u,t-1} = 0}   &\forall t\in \mathcal{T}^{\prime\prime}\\
    & \IEEEeqnarraymulticol{2}{l}{\mu^5_{u,1} - \mu^6_{u,1} + \rho_{u} \lambda_{u,1} = 0}   & \\
    & \IEEEeqnarraymulticol{2}{l}{\mu^1_{u,t} \perp \left ( x^\text{C}_t - \underline{x}^\text{C}_t \right ) = 0} & \forall t\in \mathcal{T} \label{eq:llcs1}\\
    & \IEEEeqnarraymulticol{2}{l}{\mu^2_{u,t} \perp \left ( \overline{x}^\text{C}_t - x^\text{C}_t \right ) = 0} & \forall t\in \mathcal{T} \\
    & \IEEEeqnarraymulticol{2}{l}{\mu^3_{u,t} \perp \left ( x^\text{D}_t - \underline{x}^\text{D}_t \right ) = 0} & \forall t\in \mathcal{T} \\
    & \IEEEeqnarraymulticol{2}{l}{\mu^4_{u,t} \perp \left ( \overline{x}^\text{D}_t - x^\text{D}_t \right ) = 0} & \forall t\in \mathcal{T} \\
    & \IEEEeqnarraymulticol{2}{l}{\mu^5_{u,t} \perp \left ( x^\text{S}_t - \underline{x}^\text{S}_t \right ) = 0} & \forall t\in \mathcal{T} \\
    & \IEEEeqnarraymulticol{2}{l}{\mu^6_{u,t} \perp \left ( \overline{x}^\text{S}_t - x^\text{S}_t \right ) = 0} & \forall t\in \mathcal{T} \label{eq:llcs6}\\
    & \IEEEeqnarraymulticol{2}{l}{\mu^1_{u,t}, \mu^2_{u,t}, \mu^3_{u,t},\mu^4_{u,t}, \mu^5_{u,t}, \mu^6_{u,t} \in \reals^+}&\forall t\in \mathcal{T}\\
    & \IEEEeqnarraymulticol{2}{l}{\lambda_{u,t} \in \reals}&\forall t\in \mathcal{T}
\end{IEEEeqnarray}}
\end{subequations}

A standard approach to deal with the complementary slackness conditions is a mixed-integer reformulation of \eqref{eq:llcs1} to \eqref{eq:llcs6} that can be directly implemented using off-the-shelf solvers.
The approach replaces the complementarity conditions with disjunctive constraints~\cite{FortunyAmat.1981}.
For \eqref{eq:llcs1}, the corresponding reformulation is written as:
\begin{subequations}
\begin{IEEEeqnarray}{cccr}
    \mu^1_{u,t}\; &\leq&\; z_{u,t} M & \qquad \forall t\in \mathcal{T}\,,\\
x^\text{C}_t - \underline{x}^\text{C}_t \; &\leq& \;(1-z_{u,t})M  & \qquad \forall t\in \mathcal{T}\,,
\end{IEEEeqnarray}
\end{subequations}
where $z_{u,t} \in \{0,1\}$ is a binary decision variable, and $M$ is a sufficiently large positive number. Since the underlying problem is linear, the bilevel problem can be reformulated as a mixed-integer linear problem.

\if False
\subsection{\PH{Primal-dual-equality}}
One popular approach to recast the bilevel problem as a single-level problem is to use the strong duality theorem, establishing that the optimal values of the primal and dual problems are equal for a linear programming (LP) problem.
We write the dual LP of the inner problem for an aggregated vehicle $u\in U$ as
\begin{subequations}
\begin{IEEEeqnarray}{Clll}
\maxn_{\Omega_\text{IL}}  & \IEEEeqnarraymulticol{3}{l}{\sum_{t\in \mathcal{T}} \Big [ \underline{\mu}^\text{C}_{u,t}\,\underline{x}^{\text{C}}_{u,t} - \overline{\mu}^\text{C}_{u,t}\,\overline{x}^{\text{C}}_{u,t} + \underline{\mu}^\text{D}_{u,t}\,\underline{x}^{\text{D}}_{u,t} - \overline{\mu}^\text{D}_{u,t}\,\overline{x}^{\text{D}}_{u,t}}\nonumber\\
& \IEEEeqnarraymulticol{3}{l}{\qquad + \underline{\mu}^\text{S}_{u,t}\,\underline{x}^\text{S}_{u,t} - \overline{\mu}^\text{S}_{u,t}\,\overline{x}^\text{S}_{u,t} \Big ] + \sum_{t\in \mathcal{T}^\prime} \lambda_{u,t}\,\left (\Phi^\text{DR}_{u,t} + \Phi^\text{TH}_{u,t} \right) }\nonumber\\
& \IEEEeqnarraymulticol{3}{l}{\qquad + \lambda_{u,T}\,\left (\Phi^\text{DR}_{u,T} + \Phi^\text{TH}_{u,T} + X^{S}_{u,T+1} \right)}\\
\st  & \delta^\text{C}_{u,t} + \eta^\text{C}_u\,\lambda_{u,t}  \le \Pi_t & \;\colon\; x^{\text{C}}_{u,t} & \quad \forall t\in \mathcal{T}\\
& \delta^\text{D}_{u,t} - \frac{1}{\eta^\text{D}_u}\,\lambda_{u,t}  \le -\Pi_t  & \;\colon\; x^{\text{D}}_{u,t} &\quad \forall t\in \mathcal{T}\\
& \delta^\text{S}_{u,t} + \rho_{u}\,\lambda_{u,t} - \lambda_{u,t-1}  \le 0 & \;\colon\; x^\text{S}_{u,t} & \quad\forall t\in \mathcal{T}^{\prime\prime}\\
& \delta^\text{S}_{u,1} + \rho_{u}\lambda_{u,1} \le 0 & \;\colon\; x^\text{S}_{u,1} &\\
& \underline{\mu}^\text{C}_{u,t},\underline{\mu}^\text{D}_{u,t},\underline{\mu}^\text{S}_{u,t} \in \reals^+ &&\quad\forall t\in \mathcal{T}\\
& \overline{\mu}^\text{C}_{u,t},\overline{\mu}^\text{D}_{u,t},\overline{\mu}^\text{S}_{u,t} \in \reals^+ &&\quad\forall t\in \mathcal{T}\\
& \lambda_{u,t} \in \reals && \quad\forall t\in \mathcal{T}
\end{IEEEeqnarray}
\end{subequations}
where three helper expressions are defined $\forall t \in \mathcal{T}$ as $\delta^\text{C}_{u,t} = \underline{\mu}^\text{C}_{u,t} - \overline{\mu}^\text{C}_{u,t}$, $\delta^\text{D}_{u,t} = \underline{\mu}^\text{D}_{u,t} - \overline{\mu}^\text{D}_{u,t}$, and $\delta^\text{S}_{u,t} = \underline{\mu}^\text{S}_{u,t} - \overline{\mu}^\text{S}_{u,t}$.
\par
\fi


\begin{figure}
\centering
\includegraphics[width=\columnwidth]{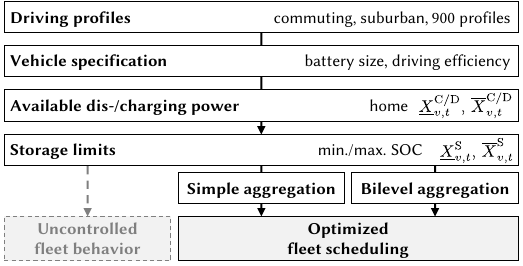}
\caption{Profile generation based on driving behavior, available charging and discharging power, and state-of-charge (SOC) limits, resulting in two EV fleet strategies: uncontrolled behavior and optimized fleet scheduling.\vspace{-1em}} 
\label{fig:gen_profiles}
\end{figure}

\section{Numerical case study}
\subsection{Case study description}
The numerical case study is designed to benchmark the proposed bilevel optimization model for AEVs against a Simple Aggregation (SA) method, a virtual storage approach.
The study evaluates the accuracy of the bilevel model by comparing it with the virtual storage-based method, using scaling factor mappings of 1h, 2h, 4h, 6h, and 24h.
These scaling factors reflect different temporal resolutions for flexibility representation.
Moreover, according to \eqref{eq:mapping}, we define the scaling factors for each weekday. 
The case study examines a three-week period in hourly resolution, resembling January 2012 in Germany.
We use projected electricity prices for 2035 to simulate future market conditions.
\par
For simplicity and to focus on the core aspects of the proposed model, the vehicle-to-grid (V2G) capability is intentionally omitted in the numerical case study and results section.
Note that we investigate an outer-level objective that only includes charging power deviations, i.e., without penalizing state-of-charge deviations.
Moreover, we use $\underline{X}^{\text{S}}_{v,t}$ instead of $\overline{X}^{\text{S}}_{v,t}$ in \eqref{eq:socminscaling}.
Again, to demonstrate the method, our dataset consists exclusively of commuters, avoiding the need for clustering and reducing heterogeneity in EV behavior.

\subsection{Individual EV profiles}
Fig.~\ref{fig:gen_profiles} illustrates the modeling framework used to generate EV demand and flexibility profiles with hourly resolution.
These can be assigned to specific locations and aggregated for representative analyses.
We use historical mobility survey data \cite{infas.2018} and traffic count data \cite{bast.2019} to create driving profiles, serving as input to determine the corresponding demand profiles and estimate flexibility potentials, i.e., available charging and discharging patterns, and EV storage limits.
These profiles imply assumptions regarding user driving behavior, the availability and power of charging infrastructure, and the primary locations where charging and discharging occur.
\par
Uncontrolled EV fleet behavior is directly derived from power availability profiles.
It considers two plug-in options: the simplest being direct charging upon arrival and a more realistic variant, low-SOC charging.
In the latter, the EV charges only when the SOC falls below a defined ``anxiety range,'' representing the SOC threshold at which users begin to feel uncomfortable and are more likely to seek the next charging opportunity.
\par
We calculate charging availability and SOC limits for each EV and then aggregate these profiles to create aggregated flexibility profiles that represent the entire EV fleet, and we consider the following cases:
\begin{itemize}
    \item \textit{Reference}: Individual EVs optimize their charging decisions as a response to the electricity price and staying close to the uncontrolled charging behavior.
    \item \textit{Simple Aggregation (SA)}: Individual EV demand and flexibility profiles are aggregated using a virtual storage approach with heuristic assumptions to avoid overestimating EV fleet flexibility. 
    \item \textit{Bilevel Aggregation}: Individual EV demand and flexibility profiles are aggregated using the proposed bilevel optimization approach with different scaling factor mappings, resulting in five AEV configurations (1h, 2h, 4h, 6h, 24h).
\end{itemize}
Note that our bilevel approach uses the reference EV unit information to identify the optimal AEV configuration, and can work with different reference cases. 

\begin{table}[htbp]
\caption{Best incumbent objective values, best lower bounds, relative gaps, and resulting charging power RMSE (w.r.t. reference of individual EVs) for different scaling factor mappings and the simple aggregation (SA).}
\label{tab:objective_values}
\centering
\begin{tabular}{|l|r|r|r|r|}
\hline
\textbf{Approach}                & \textbf{Best inc.}   & \textbf{Best bound}          & \textbf{Rel. gap}    & \textbf{RMSE in MW}\\ \hline
AEV 1h                              & 2891.69                           & 1957.95             & 32.29\,\%       & 2.926               \\
AEV 2h                              & 3734.22                           & 2921.89             & 21.75\,\%       & 3.125              \\
AEV 4h                              & 4130.89                           & 3709.89             & 10.19\,\%       & 3.164               \\
AEV 6h                              & 4300.00                           & 4162.43             & 3.20\,\%        & 3.256              \\ 
AEV 24h                             & 5666.59                           & 5609.96             & 1.00\,\%        & 3.827              \\ 
\hline
SA                     & n/a                               & n/a                 & n/a         & 18.084     \\ 
\hline
\end{tabular}
\end{table}

\begin{figure*}[htbp]
\centering
    \includegraphics[width=\textwidth]{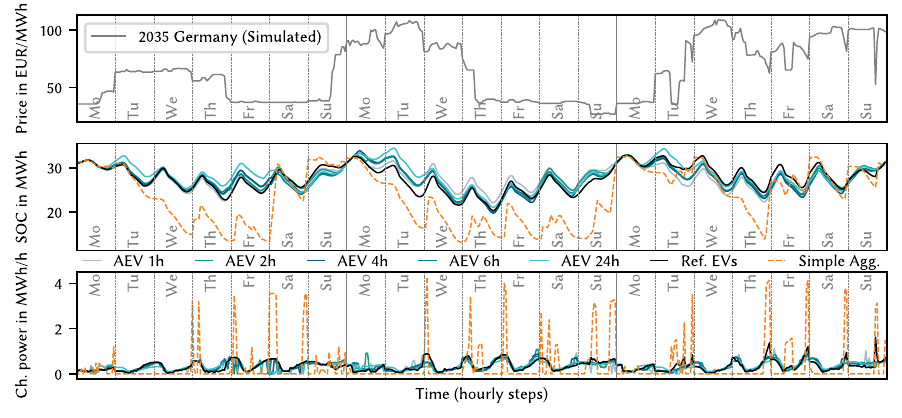}
    \caption{Electricity price, resulting Aggregated Electric Vehicle (AEV) and Simple Aggregation (SA) profiles compared to the reference EV profiles for state-of-charge (SOC) and charging power scheduling decisions during the considered three-week planning horizon.\vspace{-1em}}
    \label{fig:ev_agg_plot_price_soc_charge}
\end{figure*}

\subsection{Results and discussion}
\begin{figure}[htbp]
    \centering
    \includegraphics[width=\columnwidth]{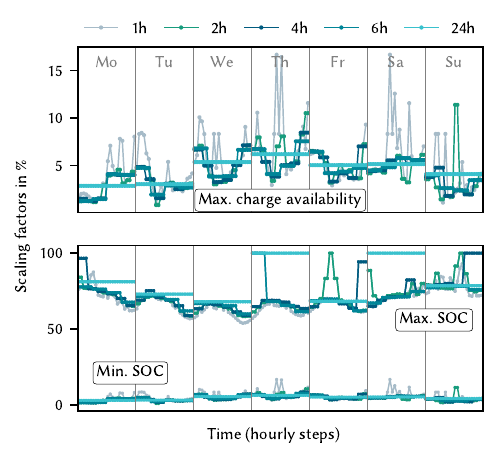}\vspace{-1em}
    \caption{Comparison of optimized scaling factors for maximum charging power availability and minimum/maximum SOC profiles for all considered scaling factor mappings and each weekday. Note that in order to obtain the final aggregated profiles shown in Figs.~\ref{fig:charge_individual} and \ref{fig:soc_individual}, the factors need to be applied to the corresponding aggregated profiles $\overline{X}^{\text{C}}_{v,t}$, $\overline{X}^{\text{S}}_{v,t}$ and $\underline{X}^{\text{S}}_{v,t}$.\vspace{-1.5em}}
    \label{fig:ev_agg_plot_scaling}
\end{figure}

Fig.~\ref{fig:ev_agg_plot_price_soc_charge} provides an overview of the key results of the numerical case study.
The electricity price exhibits several high- and low-price periods throughout the planning horizon.
Across the different scaling factor mappings, the AEV unit schedules demonstrate close alignment with the reference schedules of individual EV units, both in terms of charging power and SOC.
In contrast, the SA method shows more significant deviations from individual EV units' SOC and charging power schedules.
Due to its smaller flexibility limits and higher charging power availability, the SA method produces more pronounced charging spikes during low-price periods.
\par
Table~\ref{tab:objective_values} gives a summary of the computational results including the root-mean-square error (RMSE) for the considered time period.
While the different mappings can only be solved to optimality for one instance (relative MIP gap of 1\,\%), the different mapping variants show better objective values, i.e., smaller deviations, for the finer mappings.
Note that the best incumbents and bounds are scaled to improve the numerical stability of the underlying problem instances.
When compared to the SA approach, we find RMSE reductions of up to 78\,\% supporting the deviations observed in Fig.~\ref{fig:ev_agg_plot_price_soc_charge}. 
\par
Fig.~\ref{fig:ev_agg_plot_scaling} illustrates the resulting scaling factors for the different mappings. The restrictions on flexibility for charging availability and SOC trajectories exhibit similar patterns, with some day-to-day variations. During the three-week period under consideration, reductions in aggregate maximum charging power availability range from 0\,\% to 18\,\%. The maximum SOC profiles of the aggregated EV fleet are reduced to as low as 60\,\%, while the minimum SOC profiles range from 0\,\% to 15\,\% of the aggregate maximum SOC.
As expected, finer mappings (e.g., 1h and 2h) more accurately capture the flexibility, resulting in smaller deviations.

\section{Conclusion}
This paper develops a bilevel optimization model to enhance the flexibility aggregations of electric vehicle fleets.
The outer level minimizes scheduling deviations between the aggregated and reference EV units, while the inner level maximizes the aggregated unit's profits.
The approach is compared to a simpler aggregation method with a virtual battery approach.
The numerical case study demonstrates the bilevel model's ability to produce flexibility aggregations of electric vehicle fleets that achieve a 78\,\% reduction in charging power RMSE compared to the simple aggregation method.
The results suggest that the bilevel model can provide more accurate flexibility representations.
Moreover, the bilevel structure can be explicitly used to analyze the role of flexibility aggregators, incentivizing a desired behavior across the EV fleet.
\par
Future work should perform further case studies and sensitivity analyses, including the vehicle-to-grid options, different user types, various price scenarios, fast charging, charging at work, and the impact of day-specific scaling factors to take into account the fact that the driving patterns and thus the charging behavior differ between the days of the week.


%

\appendices
\section{Reference and aggregated results}
\begin{figure*}
\centering
\includegraphics[width=0.99\textwidth]{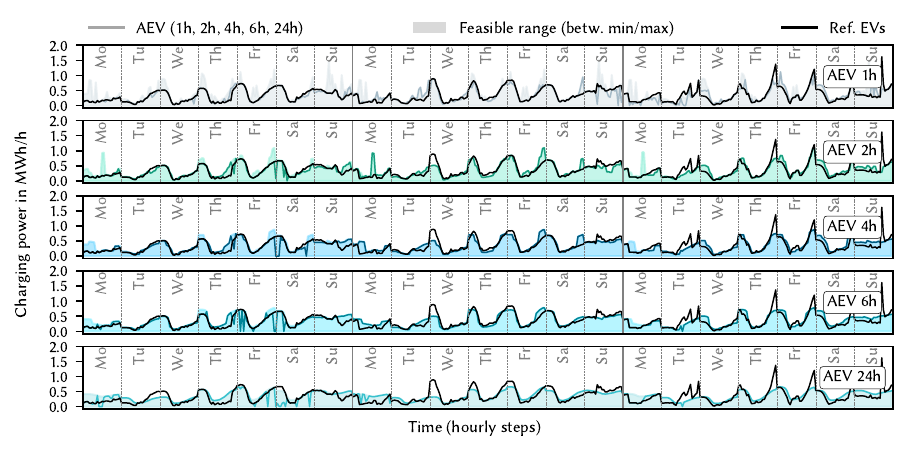}
\caption{Reference (black line) and resulting AEV charging profiles and feasibility space for each scaling factor mapping (different colors for 1h, 2h, 4h, 6h, 24h), own illustration based on own computations.}
\label{fig:charge_individual}
\end{figure*}

\begin{figure*}
\centering
\includegraphics[width=0.99\textwidth]{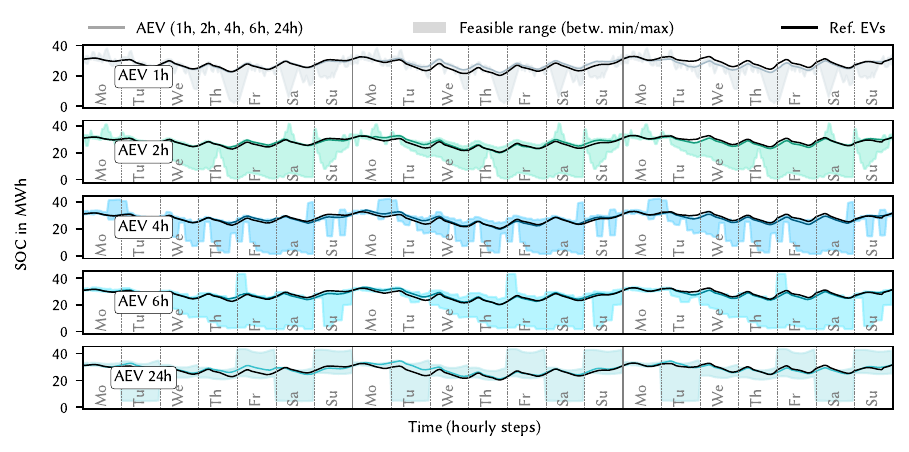}
\caption{Reference (black line) and resulting AEV state-of-charge (SOC) profiles and feasibility space for each scaling factor mapping (different colors for 1h, 2h, 4h, 6h, 24h), own illustration based on own computations.}
\label{fig:soc_individual}
\end{figure*}

Figs.~\ref{fig:charge_individual} and \ref{fig:soc_individual} illustrate the reference and aggregate profiles and feasibility space for each scaling factor mapping (1h, 2h, 4h, 6h, 24h) for the charging power and SOC trajectories, respectively.


\section*{Acknowledgment}
The authors gratefully acknowledge Gurobi Optimization LLC for their \textit{Gurobi Optimizer}, which was instrumental in solving the optimization problems presented in this study.

\ifCLASSOPTIONcaptionsoff
  \newpage
\fi

\bibliographystyle{IEEEtran}

\bibliography{mybibfile}

\begin{thebibliography}{10}
\providecommand{\url}[1]{#1}
\csname url@samestyle\endcsname
\providecommand{\newblock}{\relax}
\providecommand{\bibinfo}[2]{#2}
\providecommand{\BIBentrySTDinterwordspacing}{\spaceskip=0pt\relax}
\providecommand{\BIBentryALTinterwordstretchfactor}{4}
\providecommand{\BIBentryALTinterwordspacing}{\spaceskip=\fontdimen2\font plus
\BIBentryALTinterwordstretchfactor\fontdimen3\font minus
  \fontdimen4\font\relax}
\providecommand{\BIBforeignlanguage}[2]{{%
\expandafter\ifx\csname l@#1\endcsname\relax
\typeout{** WARNING: IEEEtran.bst: No hyphenation pattern has been}%
\typeout{** loaded for the language `#1'. Using the pattern for}%
\typeout{** the default language instead.}%
\else
\language=\csname l@#1\endcsname
\fi
#2}}
\providecommand{\BIBdecl}{\relax}
\BIBdecl

\bibitem{Braun.2024}
M.~Braun, C.~Gruhl, C.~A. Hans, P.~H{\"a}rtel, C.~Scholz, B.~Sick, M.~Siefert,
  F.~Steinke, O.~Stursberg, and S.~{Wende-von Berg}, ``Predictions and decision
  making for resilient intelligent sustainable energy systems,'' in \emph{2024
  IEEE PES Innovative Smart Grid Technologies Europe (ISGT EUROPE)}.\hskip 1em
  plus 0.5em minus 0.4em\relax IEEE, 2024, pp. 1--5.

\bibitem{Rogers.1991}
D.~F. Rogers, R.~D. Plante, R.~T. Wong, and J.~R. Evans, ``Aggregation and
  disaggregation techniques and methodology in optimization,'' \emph{Operations
  Research}, vol.~39, no.~4, pp. 553--582, 1991.

\bibitem{Hartel.2021}
\BIBentryALTinterwordspacing
P.~H{\"a}rtel, \emph{Offshore Grids in Low-Carbon Energy Systems: Long-Term
  Transmission Expansion Planning in Energy Systems with Cross-Sectoral
  Integration using Decomposition Algorithms and Aggregation Methods for
  Large-Scale Optimisation Problems}.\hskip 1em plus 0.5em minus 0.4em\relax
  Stuttgart: {Fraunhofer Verlag}, 2021. [Online]. Available:
  \url{https://doi.org/10.24406/publica-fhg-283522}
\BIBentrySTDinterwordspacing

\bibitem{Muessel.2023}
J.~Muessel, O.~Ruhnau, and R.~Madlener, ``Accurate and scalable representation
  of electric vehicles in energy system models: A virtual storage-based
  aggregation approach,'' \emph{iScience}, vol.~26, no.~10, p. 107816, 2023.

\bibitem{Trost.2017}
\BIBentryALTinterwordspacing
T.~Trost, \emph{{Erneuerbare Mobilit{\"a}t im motorisierten Individualverkehr:
  Modellgest{\"u}tzte Szenarioanalyse der Marktdiffusion alternativer
  Fahrzeugantriebe und deren Auswirkungen auf das
  Energieversorgungssystem}}.\hskip 1em plus 0.5em minus 0.4em\relax Stuttgart:
  {Fraunhofer Verlag}, 2017. [Online]. Available:
  \url{https://doi.org/10.24406/publica-fhg-281311}
\BIBentrySTDinterwordspacing

\bibitem{Pertl.2019}
M.~Pertl, F.~Carducci, M.~Tabone, M.~Marinelli, S.~Kiliccote, and E.~C. Kara,
  ``An equivalent time-variant storage model to harness ev flexibility:
  Forecast and aggregation,'' \emph{IEEE Transactions on Industrial
  Informatics}, vol.~15, no.~4, pp. 1899--1910, 2019.

\bibitem{Geng.2023}
\BIBentryALTinterwordspacing
S.~Geng, T.~Lee, D.~Mallapragada, and A.~Botterud, ``An integer clustering
  approach for modeling large-scale ev fleets with guaranteed performance.''
  [Online]. Available: \url{https://arxiv.org/pdf/2310.02208.pdf}
\BIBentrySTDinterwordspacing

\bibitem{Mukhi.2024}
K.~Mukhi, C.~Qu, P.~You, and A.~Abate, ``Robust aggregation of electric vehicle
  flexiblity.''

\bibitem{Almeida.2013}
K.~C. Almeida and A.~J. Conejo, ``Medium-term power dispatch in predominantly
  hydro systems: An equilibrium approach,'' \emph{IEEE Transactions on Power
  Systems}, vol.~28, no.~3, pp. 2384--2394, 2013.

\bibitem{Risberg.2017}
D.~Risberg and L.~Söder, ``Hydro power equivalents of complex river systems,''
  in \emph{2017 IEEE Manchester PowerTech}, 2017, pp. 1--6.

\bibitem{Blom.2020}
E.~Blom, L.~S{\"o}der, and D.~Risberg, ``Performance of multi-scenario
  equivalent hydropower models,'' \emph{Electric Power Systems Research}, vol.
  187, p. 106486, 2020.

\bibitem{Schmitz.2023}
R.~Schmitz, C.~{\O}. Naversen, and P.~H{\"a}rtel, ``Influence of hydrogen
  import prices on hydropower systems in climate-neutral {E}urope,''
  \emph{Energy Systems}, 2023.

\bibitem{Dempe.2015}
S.~Dempe, V.~Kalashnikov, G.~A. P{\'e}rez-Vald{\'e}s, and N.~Kalashnykova,
  \emph{Bilevel Programming Problems: Theory, Algorithms and Applications to
  Energy Networks}, ser. Energy systems.\hskip 1em plus 0.5em minus 0.4em\relax
  Berlin, Heidelberg and s.l.: {Springer Berlin Heidelberg}, 2015.

\bibitem{FortunyAmat.1981}
J.~Fortuny-Amat and B.~McCarl, ``A representation and economic interpretation
  of a two-level programming problem,'' \emph{The Journal of the Operational
  Research Society}, vol.~32, no.~9, pp. 783--792, 1981.

\bibitem{infas.2018}
{Infas Institut f{\"u}r Sozialwissenschaft}, ``Mobilit{\"a}t in
  {D}eutschland,'' \url{http://www.mobilitaet-in-deutschland.de}, Online, 2018.

\bibitem{bast.2019}
{Bundesanstalt f{\"u}r Stra{\ss}enwesen}, ``Automatische {Z}{\"a}hlstellen auf
  {A}utobahnen und {B}undesstra{\ss}en,''
  \url{https://www.bast.de/BASt\_2017/DE/Verkehrstechnik/Fachthemen/v2-verkehrszaehlung/zaehl\_node.html},
  2019, accessed: 26.11.2021.

\end{thebibliography}

\end{document}